\documentclass[12pt,a4paper]{article}
\usepackage{german, latexsym, a4,diagrams,colordvi}

\author{Projektarbeit von Martin Kohls\\
Betreuer: Prof. Dr. Gregor Kemper}
\title{Konstruktion von Invariantenringen ohne die Cohen-Macaulay Eigenschaft}
\date{17. Oktober 2003}
\begin{document}
\maketitle
\begin{abstract}
\noindent In dieser Projektausarbeitung werden f"ur einen algebraisch abgeschlos\-senen K"orper $K$ positiver Charakteristik und die Gruppen $G=SL_{n}(K)$ bzw. $GL_{n}(K)$ Beispiele von $G$-Moduln $V$ konstruiert, so dass der Invariantenring $K[V]^{G}$ nicht Cohen-Macaulay ist.
\end{abstract}

\tableofcontents

\section*{Einleitung}
Die Cohen-Macaulay Eigenschaft spielt eine wichtige Rolle in der Invarianten- theorie. Der ber"uhmte Satz von Hochster und Roberts \cite{HoRo} besagt z.B., dass f"ur eine linear reduktive Gruppe $G$ und einen $G$-Modul $V$ der Invariantenring $K[V]^{G}$ Cohen-Macaulay ist. In der Arbeit \cite[S. 1026]{HoEa} dagegen bemerken Hochster und Eagon:

\textit{We know of no example in which $G$ is reductive or connected [...] and $R^{G}$ is not Cohen-Macaulay.}\\
In dieser Arbeit werden solche Beispiele f"ur $R=K[V]$ konstruiert.

Nach Kemper \cite{Kemper1} gilt auch die folgende Umkehrung des Resultats von Hochster und Roberts: Wenn $G$ reduktiv und $K[V]^{G}$ f"ur jeden $G$-Modul $V$ Cohen-Macaulay ist, dann ist $G$ sogar linear reduktiv. Den dort gegebenen Widerspruchsbeweis verwenden wir f"ur die Konstruktion eines nicht Cohen-Macaulay Invarianten\-rings $K[V]^{G}$.
Da in Charakteristik 0 die Eigenschaften reduktiv und linear re\-duk\-tiv "aquivalent sind, k"onnen wir dabei aufgrund des Resultats von Hochster und Roberts  positive Charakteristik voraussetzen.

Im Folgenden bezeichnen wir daher mit $K$ stets einen algebraisch abge- schlossenen K"orper mit positiver Charakteristik $\textrm{char } K=p$ und mit $G$ stets eine lineare algebraische Gruppe. Wir geben zun"achst den Teil der Arbeit \cite{Kemper1} wieder, der zu folgendem Resultat f"uhrt:\\
\textit{Ist $G$ eine reduktive Gruppe, $0 \rightarrow U  \rightarrow  \tilde{U} \rightarrow K \rightarrow 0$ eine kurze exakte Sequenz von $G$-Moduln, die nicht zerf"allt, so ist mit $V:=U^{*}\oplus\tilde{U}\oplus\tilde{U}\oplus\tilde{U}$ der Invariantenring $K[V]^{G}$ nicht Cohen-Macaulay.\\}
Das eigentliche Resultat des Projekts ist die Angabe solcher Sequenzen, womit die Arbeit schlie\ss t.

\section{Erste Kohomologie algebraischer Gruppen}
Sei $V$ ein $G$-Modul. Ein Morphismus von affinen Variet"aten\footnote{d.h. eine durch Polynome in den Koeffizienten eines Gruppenelements $\sigma \in G \subseteq K^{r}$ gegebene Abbildung} $g: G \rightarrow V, \quad \sigma \mapsto g_{\sigma}$ hei\ss t \textit{1-Kozyklus}, falls $g_{\sigma\tau}=\sigma(g_{\tau})+g_{\sigma}$ f"ur alle $\sigma,\tau \in G$. Die additive Gruppe aller 1-Kozyklen (die zugleich ein $K$-Vektorraum ist) wird mit $Z^{1}(G,V)$ bezeichnet. F"ur ein $v \in V$ ist durch $\sigma \mapsto (\sigma - 1)v:=\sigma(v)-v$ ein spezieller 1-Kozyklus gegeben. Die Untergruppe dieser 1-Kozyklen wird mit $B^{1}(G,V)$ bezeichnet, und die zugeh"orige Faktorgruppe $Z^{1}(G,V)/B^{1}(G,V)$ mit $H^{1}(G,V)$.\\

\noindent {\bf Bemerkung 1.} \textit{Sei
\[
0 \rightarrow V  \hookrightarrow  \tilde{V}  \stackrel{\pi}{\rightarrow} K \rightarrow 0
\]
eine kurze exakte Sequenz von $G$-Moduln. \footnote{$K$ soll dabei stets triviale $G$-Operation haben} Diese zerf"allt genau dann, wenn f"ur ein (und dann alle) $v_{0} \in \pi^{-1}(1)$ und den durch $g_{\sigma}:=(\sigma-1)v_{0}$ definierten Kozyklus $g \in Z^{1}(G,V)$ gilt, dass $g$ sogar in $B^{1}(G,V)$ liegt.}\\

\noindent \textit{Beweis.} Wenn die Sequenz zerf"allt, so hat $V$ ein Komplement $W$ in $\tilde{V}$, das wegen $\textrm{dim }\tilde{V} = \textrm{dim Im }\pi + \textrm{dim Kern }\pi = 1 + \textrm{dim  }V$ eindimensional ist. Ist $v_{0}=v+w$ mit $v \in V, w \in W$, so gibt es also zu $\sigma \in G$ ein $\lambda \in K$ mit $\sigma w=\lambda w$. Aus $1=\pi(v_{0})=\pi(w)=\sigma \pi(w)=\pi(\sigma w)=\pi(\lambda w)=\lambda \pi(w)$ folgt $\lambda=1$, und $w$ ist $G$-invariant. Es folgt $g_{\sigma}=(\sigma -1)v$, also $g \in B^{1}(G,V)$.\\

Ist umgekehrt $g \in B^{1}(G,V)$, also $g_{\sigma}=(\sigma-1)v$ mit $v \in V$, so ist $K(v_{0}-v)$ ein $G$-invariantes Komplement zu $V$. \hfill $\Box$\\

\noindent Es ist n"utzlich, sich mit der Darstellung von $G$ auf $\tilde{V}$ vertraut zu machen: Wir erg"anzen $v_{0}$ mit Hilfe einer Basis von $V$ zu einer Basis von $\tilde{V}$. Ist dann $A_{\sigma}$ die Darstellung von $\sigma$ auf $V$ und identifizieren wir $g_{\sigma}$ mit seinem Koordinatenvektor bzgl. der Basis von $V$, so hat $\sigma$ bzgl. der Basis von $\tilde{V}$ die Darstellung  
$\left( \begin{array}{cc}
A_{\sigma} & g_{\sigma} \\
0 & 1
\end{array} \right)$. Insbesondere finden wir in der Identit"at
\[
\left( \begin{array}{cc}
A_{\sigma\tau} & g_{\sigma\tau} \\
0 & 1
\end{array} \right)
=
\left( \begin{array}{cc}
A_{\sigma} & g_{\sigma} \\
0 & 1
\end{array} \right)
\cdot
\left( \begin{array}{cc}
A_{\tau} & g_{\tau} \\
0 & 1
\end{array} \right)
=
\left( \begin{array}{cc}
A_{\sigma}A_{\tau} & A_{\sigma}g_{\tau}+g_{\sigma} \\
0 & 1
\end{array} \right)
\]
die Kozyklus-Eigenschaft von $g$ wieder. Daher kann man aus gegebenem $g \in Z^{1}(G,V)$ durch $\tilde{V}:=V \oplus K$ und  $\sigma (v,\lambda):=(\sigma v+\lambda g_{\sigma},\lambda) \cong 
\left( \begin{array}{cc}
A_{\sigma} & g_{\sigma} \\
0 & 1
\end{array} \right)
\left( \begin{array}{c}
v\\ \lambda
\end{array} \right)
$
sowie $\pi(v,\lambda):=\lambda$ eine kurze exakte Sequenz definieren. Wir sagen, $\tilde{V}$ wird von $g$ induziert. Ist $(v,1) \in \pi^{-1}(1)$, so ist $(\sigma-1)(v,1)=(g_{\sigma}+\sigma (v) -v,0)$, so dass man aus dieser Sequenz die Restklasse von $g$ in $H^{1}(G,V)$ zur"uckgewinnen kann.\\

\noindent Wir bezeichnen im Folgenden die durch ein $\sigma \in G$ durch Operation induzierte lineare Abbildung ebenfalls mit $\sigma$, auch dann wenn damit verschiedene lineare Abbildungen auf zwei $G$-Moduln $V,W$ mit dem selben Symbol bezeichnet werden. Mit dieser Notation wird $\textrm{Hom}_{K}(V,W)$ zu einem $G$-Modul durch die Operation 
\[
G \times \textrm{Hom}_{K}(V,W) \rightarrow \textrm{Hom}_{K}(V,W), \quad (\sigma,f) \mapsto \sigma f := \sigma \circ f \circ \sigma^{-1}.
\]
Da $G$ auf $K$ stets trivial operieren soll, haben wir damit insbesondere f"ur $\varphi \in V^{*}=\textrm{Hom}_{K}(V,K)$ die Operation $\sigma \varphi = \varphi \circ \sigma^{-1}$.\\
Es ist
\[
f \in \textrm{Hom}_{K}(V,W)^{G} \Leftrightarrow \sigma f=f \textrm{ bzw. } f  \circ \sigma = \sigma \circ f  \quad \forall  \sigma \in G,
\]
also gilt $\textrm{Hom}_{K}(V,W)^{G}=\textrm{Hom}_{G}(V,W)$.\\

\noindent F"ur sp"ater ben"otigen wir noch zwei kanonische Homomorphismen. Aufgrund der endlichen Dimension ist 
\[
f: V \rightarrow V^{**}, \quad v \mapsto f_{v}: \quad V^{*} \ni \varphi \mapsto f_{v}(\varphi ):=\varphi (v)
\]
ein $K$-Isomorphismus. Wegen
\[
f_{\sigma v}(\varphi )=\varphi (\sigma v)=(\sigma ^{-1} \varphi )(v)=f_{v}(\sigma^{-1}\varphi )=(\sigma f_{v})(\varphi )
\]
gilt $f_{\sigma v}=\sigma f_{v}$ und $f$ ist auch ein $G$-Homomorphismus. Es gilt also $V\cong V^{**}$ auch als $G$-Moduln.\\

\noindent Durch lineare Fortsetzung wird durch
\[
W\otimes V \rightarrow \textrm{Hom}_{K}(V^{*},W), \qquad w \otimes v \mapsto f_{w \otimes v}: \quad V^{*} \ni \varphi \mapsto \varphi(v)w
\] 
ein weiterer kanonischer $K$-Isomorphismus gegeben, der wegen
\[
\begin{array}{c}
f_{\sigma w \otimes \sigma v}(\varphi)=
 \varphi(\sigma v)\sigma w=\sigma( (\sigma^{-1}\varphi)(v)w)=\\
\sigma(f_{w \otimes v}(\sigma^{-1}\varphi))= 
(\sigma \circ f_{w \otimes v} \circ \sigma^{-1})(\varphi)
=(\sigma f_{w \otimes v})(\varphi)
\end{array}
\]
die Gleichung $f_{\sigma (w \otimes v)}=\sigma f_{w \otimes v}$ erf"ullt, also ein $G$-Homomorphismus ist.\\

\noindent {\bf Proposition 2. }\textit{$G$ ist genau dann linear reduktiv, wenn $H^{1}(G,V)=0$ f"ur jeden $G$-Modul $V$ gilt.}\\

\noindent \textit{Beweis. } Wenn $G$ linear reduktiv ist, zerf"allt jede kurze exakte Sequenz
$0 \rightarrow V  \rightarrow  \tilde{V} \rightarrow K \rightarrow 0$, also gilt nach Bemerkung 1, dass $H^{1}(G,V)=0$.

Sei umgekehrt $H^{1}(G,V)=0$ f"ur jeden $G$-Modul $V$. Dann ist jede kurze exakte Sequenz wie oben zerfallend, und f"ur jeden Epimorphismus $F:\quad V \rightarrow W$ von $G$-Moduln ist auch die Restriktion $V^{G} \rightarrow W^{G}$ surjektiv: Denn f"ur $0 \ne w \in W^{G}$ ist mit $V_{w}=F^{-1}(<w>)$ eine kurze exakte Sequenz $0 \rightarrow \textrm{Ker }F  \rightarrow  V_{w} \rightarrow <w>\cong K \rightarrow 0$ gegeben, die nach Voraussetzung zerf"allt und daher ein Urbild von w in $V^{G}$ liefert. Sei nun eine beliebige kurze exakte Sequenz
\begin{equation} \label{seq}
0 \rightarrow U \stackrel{\varepsilon}{\rightarrow} V \stackrel{\pi}{\rightarrow} W \rightarrow 0
\end{equation}
von $G$-Moduln gegeben. Diese gibt mittels
\[
\Gamma: \textrm{Hom}_{K}(W,U) \rightarrow \textrm{Hom}_{K}(V,U),\quad g \mapsto g \circ \pi
\]
\[
\Lambda: \textrm{Hom}_{K}(V,U) \rightarrow \textrm{Hom}_{K}(U,U), \quad f \mapsto f \circ \varepsilon
\]
Anlass zu folgender kurzen exakten Sequenz:
\[
0 \rightarrow \textrm{Hom}_{K}(W,U) \stackrel{\Gamma}{\rightarrow} \textrm{Hom}_{K}(V,U) \stackrel{\Lambda}{\rightarrow} \textrm{Hom}_{K}(U,U) \rightarrow 0.
\]
Dabei folgt die Injektivit"at von $\Gamma$ aus der Surjektivit"at von $\pi$, und die Surjektivit"at von $\Lambda$ aus der Injektivit"at von $\varepsilon$. Aus $\pi \circ \varepsilon =0$ folgt $\textrm{Im }\Gamma \subseteq \textrm{Ker }\Lambda$. Die umgekehrte Inklusion folgt wegen Dimensionsgr"unden aus einem Liftungsargument. \footnote{Aus $f \in \textrm{Ker }\Lambda$ folgt $f \circ \varepsilon=0$, also $f|\textrm{Im }\varepsilon=f|\textrm{Ker }\pi=0.$ Da $\pi$ surjektiv, kann man $f=g \circ \pi$ sogar mit eindeutig bestimmten $g$ schreiben.}\\
Wie wir gesehen haben, ist dann jedenfalls auch die Restriktion
$\textrm{Hom}_{G}(V,U) \stackrel{\Lambda}{\rightarrow} \textrm{Hom}_{G}(U,U)$ surjektiv, insbesondere gibt es ein Urbild von $\textrm{id}_{U}$, also ein $f$ mit $f \circ \varepsilon=\textrm{id}_{U}$. Dies zeigt, dass (\ref{seq}) zerf"allt, und $U$ hat ein Komplement in $V$, d.h. $G$ ist linear reduktiv. \hfill $\Box$\\

\noindent Sei $g \in Z^{1}(G,V)$, wobei wir die zugeh"orige Restklasse von $g$ in $H^{1}(G,V)$ ebenfalls mit $g$ bezeichnen. Ist $W$ ein weiterer $G$-Modul, $w \in W^{G}$, so ist $\sigma \mapsto w \otimes g_{\sigma}$ wegen $\sigma w =w \quad \forall \sigma \in G$ aus $Z^{1}(G,W \otimes V)$. Die zugeh"orige Restklasse in $H^{1}(G, W \otimes V)$ bezeichnen wir mit $w \otimes g$. Sie ist unabh"angig von dem Repr"asentanten von $g$ in $Z^{1}(G,V)$, da
\[
w \otimes (g_{\sigma}+\sigma (v)-v)=w \otimes g_{\sigma}+\sigma(w \otimes v)- w \otimes v 
\]
wegen $w \in W^{G}$.\\

\noindent {\bf Proposition 3.} \textit{Sei $V$ ein $G$-Modul, $g \in H^{1}(G,V)$ und $\tilde{V}$ der von einem Repr"asentanten von $g$ in $Z^{1}(G,V)$  induzierte $G$-Modul. Dann gibt es einen $G$-Modul $W$ (n"amlich $W=\tilde{V}^{*}$) und ein $0 \ne w \in W^{G}$ mit $w\otimes g=0$ in $H^{1}(G, W \otimes V)$.}\\

\noindent \textit{Beweis.\footnote{siehe auch \cite{Kemper2}} } Die von $g$ induzierte kurze exakte Sequenz bezeichnen wir mit
\[
0 \rightarrow V  \rightarrow  \tilde{V} \stackrel{\pi}{\rightarrow} K \rightarrow 0.
\]
Sei zun"achst $W$ ein beliebiger $G$-Modul, und $0 \ne w \in W^{G}$. Wir leiten Bedingungen f"ur $w \otimes g=0$ her. Wir betrachten den $G$-Homomorphismus
\[
\varphi:  K \rightarrow W,\qquad c \mapsto cw
\]
und den "`Pullback"'
\[
M:=(W \otimes \tilde{V}) \times_{W} K:=\left\{ (a,c) \in (W \otimes \tilde{V})\oplus K : (\textrm{id}_{W} \otimes  \pi )(a)=\varphi (c) \right\}.
\]
Wir berechnen $M$ explizit. Sei dazu $\{w_{1},...,w_{n}\}$ Basis von $W$ mit $w_{1}=w$, und $\{v_{1},...,v_{m}\}$ Basis von $V$. Dann hat $a \in W \otimes \tilde{V}$ eine Darstellung
\[
 a=\sum_{i=1..n,  j=1..m}a_{ij} w_{i} \otimes (v_{j},0)+ \sum_{i=1..n} b_{i} w_{i} \otimes (0,1) \qquad  \textrm{mit } a_{ij},b_{i} \in K.
\]
Also haben wir die Bedingung
\[
(\textrm{id}_{W} \otimes \pi ) (a)= \sum_{i=1..n} b_{i} w_{i} \stackrel{!}{=} \varphi(c)=cw_{1}
\]
und damit $b_{i}=c\delta_{1i}$. Es gilt also
\[
(a,c) \in M \Leftrightarrow  a=\sum_{i=1..n,  j=1..m}a_{ij} w_{i} \otimes (v_{j},0)+ c w_{1} \otimes (0,1). 
\]
Da
\[
\sigma(c w_{1} \otimes (0,1))=c w_{1} \otimes (g_{\sigma},1)=c w_{1} \otimes (g_{\sigma},0)+c w_{1} \otimes(0,1),
\]
wird $M$ durch $\sigma(a,c):=(\sigma a,c)$ zu einem $G$-Modul.\\
Offenbar ist $M\cong \widetilde{W \otimes V}$, wobei $\widetilde{W \otimes V}$ aus $W \otimes V$ mittels $(w_{1} \otimes g_{\sigma}) \in Z^{1}(G, W \otimes V)$ entsteht. Dabei ist $\widetilde{W \otimes V} \cong M$ gegeben durch
\[
\left(\sum a_{ij} w_{i} \otimes v_{j},c\right) \mapsto \sum a_{ij} w_{i} \otimes (v_{j},0)+c w_{1} \otimes(0,1).
\]
Daher ist mit $\tilde{\pi}(a,c)=c$ durch
\begin{equation} \label{seq2}
0 \rightarrow W \otimes V  \rightarrow  M \stackrel{\tilde{\pi}}{\rightarrow} K \rightarrow 0
\end{equation}
eine kurze exakte Sequenz gegeben, die genau dann zerf"allt, wenn $w \otimes g=0$ gilt.\\
Andererseits zerf"allt (\ref{seq2}) genau dann, wenn es einen $G$-Homomorphismus $p':K \rightarrow W \otimes \tilde{V}$ gibt, so dass folgendes Diagramm
\begin{diagram}
	&			& K\\
	&\ldTo^{p'}		& \dTo>{\varphi} \\
W \otimes \tilde{V} & \rTo^{\textrm{id}_{W} \otimes \pi}& W\\
\end{diagram}
kommutiert. In diesem Fall gilt n"amlich mit 
\[
p:K \rightarrow M, \quad c \mapsto (p'(c),c) \stackrel{!}{\in} M
\]
(wegen $\textrm{id}_{W} \otimes \pi (p'(c))=\varphi(c)$), dass $\tilde{\pi}(p(c))=c$, also $\tilde{\pi} \circ p=\textrm{id}_{K}$, d.h. (\ref{seq2}) zerf"allt.\\
Der kanonische Homomorphismus $W \otimes \tilde{V} \rightarrow \textrm{Hom}_{K}(\tilde{V}^{*},W)$ und der $G$-Homomorphismus 
\footnote{Weil $\sigma f= \sigma \circ f \circ \sigma^{-1} \mapsto \sigma (f( \sigma^{-1}(\pi))=\sigma (f(\pi )),$ da $\pi \in \textrm{Hom}_{G}(\tilde{V},K)=\textrm{Hom}_{K}(\tilde{V},K)^{G}$. }
\[
\textrm{Hom}_{K}(\tilde{V}^{*},W) \rightarrow W, \quad f \mapsto f(\pi)
\]
lassen folgendes Diagramm kommutieren,
\begin{diagram}
W \otimes \tilde{V}	& \rTo^{\textrm{id}_{W} \otimes \pi}	& W \\
				& \rdTo<{\cong}					& \uTo\\
				&							& \textrm{Hom}_{K}(\tilde{V}^{*},W)\\
\end{diagram}
wie man am leichtesten mit Tensoren "uberpr"uft:
\begin{diagram}
w \otimes (v,\lambda)	& \rTo^{\textrm{id}_{W} \otimes \pi}	& \lambda w \\
				& \rdTo<{\cong}					& \uTo>(\psi=\pi)\\
				&							& (\psi \mapsto \psi(v,\lambda) w)\\
\end{diagram}
Durch Zusammenf"ugen beider Diagramme sieht man, dass (\ref{seq2}) genau dann zerf"allt (bzw. $w \otimes g =0$ gilt), wenn es einen $G$-Homomorphismus $K \rightarrow \textrm{Hom}_{K}(\tilde{V}^{*},W)$ gibt, so dass
\begin{diagram}
						&			&K\\
						& \ldTo		&\dTo>{\varphi}\\
\textrm{Hom}_{K}(\tilde{V}^{*},W)	& \rTo 		&W\\
\end{diagram}
kommutiert. Bezeichnet man das Bild von $1$ unter diesem Homomorphismus mit $f$, und beachtet man, dass ein solcher Homomorphismus eindeutig durch sein ($G$-invariantes) Bild von $1$ bestimmt ist, so sieht man, dass die letzte Bedingung "aquivalent ist zur Existenz eines
\[
f \in \textrm{Hom}_{K}(\tilde{V}^{*},W)^{G}=\textrm{Hom}_{G}(\tilde{V}^{*},W) \textrm{   mit } f(\pi)=\varphi(1)=w.
\]
Also k"onnen wir $W=\tilde{V}^{*},\quad f=\textrm{id}, \quad w:=f(\pi)=\pi$ w"ahlen. \hfill$\Box$

\section{Eine Eigenschaft reduktiver Gruppen}
F"ur die Konstruktion eines nicht Cohen-Macaulay Invariantenrings $K[V]^{G}$ ben"otigen wir nur ein einziges Lemma "uber reduktive Gruppen.  Bevor wir es formulieren, stellen wir einige Aussagen "uber Parametersysteme zusammen.\\
Sei dazu im folgenden $R$ eine noethersche, zusammenh"angende, graduierte, nullteilerfreie Algebra "uber $K$.\\
Dann bilden homogene Elemente $f_{1},...,f_{r} \in R$ ein \emph{partielles homogenes parametersystem (phsop)}, wenn sie sich zu einem homogenen Parameter\-sy\-stem (hsop) $f_{1},...,f_{n} \in R$ erg"anzen lassen, d.h. $f_{1},...,f_{n}$  sind homogen und algebra\-isch unabh"angig "uber $K$, und $R$ ist endlich erzeugt als Modul "uber $A=K[f_{1},...,f_{n}]$. $R$ ist \emph{Cohen-Macaulay}, falls $R$ sogar frei "uber $A$ ist (f"ur ein und dann alle hsop).\\

$\{f_{1},...f_{r}\}$ hei\ss t eine \emph{regul"are Sequenz}, falls $\forall i=1..r$ gilt:
\[
\forall g\in R: \qquad f_{i}g \in (f_{1},...,f_{i-1}) \Rightarrow g \in (f_{1},...,f_{i-1}).
\]
Jede regul"are Sequenz ist ein phsop, und wenn $R$ Cohen-Macaulay ist, so ist jedes phsop eine regul"are Sequenz.\\
Letztere Eigenschaft wird f"ur die Konstruktion eines nicht Cohen-Macaulay Invariantenrings entscheidend sein, ebenso wie folgendes Lemma "uber die ben"otigte Eigenschaft reduktiver Gruppen, welches wir ohne Beweis for\-mu\-lier\-en. (Siehe \cite{Kemper1} f"ur einen Beweis.)\\
\\
\noindent {\bf Lemma 4. } \textit{Sei $G$ eine reduktive Gruppe und  $V$ ein $G$-Modul. Wenn $a_{1},...,a_{k} \in K[V]^{G}$ ein phsop in $K[V]$ bilden, dann auch in $K[V]^{G}$.}\\

\section{Konstruktion eines nicht Cohen-Macaulay Invariantenrings}

\noindent{\bf Proposition 5. } \textit{ Sei $G$ eine reduktive Gruppe und $V$ ein $G$-Modul, so dass $K[V]^{G}$ Cohen-Macaulay ist. Falls die Elemente $a_{1},a_{2},a_{3} \in K[V]^{G}$ ein phsop in $K[V]$ bilden, dann ist die durch die Multiplikation mit $(a_{1},a_{2},a_{3})$ gegebene Abbildung
\[
H^{1}(G,K[V]) \rightarrow H^{1}(G,K[V]^{3})
\]
injektiv. }\\
\\
\noindent \textit{Beweis. } Nach Lemma 4 bilden $a_{1},a_{2},a_{3}$ auch ein phsop in dem CM Ring $K[V]^{G}$, also dort sogar eine regul"are Sequenz. Sei nun $g \in Z^{1}(G,K[V])$ ein Kozyklus, dessen Restklasse in $H^{1}(G,K[V])$ durch Multiplikation mit $(a_{1},a_{2},a_{3})$ auf $0$ abgebildet wird. \footnote{ Aus $a_{i} \in K[V]^{G}$ folgt $a_{i} g_{\sigma \tau}=a_{i} \sigma (g_{\tau})+a_{i} g_{\sigma}= \sigma (a_{i} g_{\tau})+ a_{i} g_{\sigma}$, d.h $a_{i} g \in Z^{1}(G,K[V])$.} Dann gibt es $b_{i} \in K[V]$  mit $(\sigma -1)b_{i}=a_{i} g_{\sigma} \quad \forall \sigma \in G,\quad i=1,2,3$. Sei $u_{ij}=a_{i}b_{j}-a_{j}b_{i}$ f"ur $1 \le i < j \le 3$.  Offenbar ist $u_{ij} \in K[V]^{G}$, und es gilt
\[
u_{23}a_{1}-u_{13}a_{2}+u_{12}a_{3}=
\left|
\begin{array}{ccc}
a_{1} & a_{2} & a_{3}\\
a_{1} & a_{2} & a_{3}\\
b_{1} & b_{2} & b_{3}\\
\end{array}
\right|=0.
\]
Aufgrund der Regularit"at von $a_{1},a_{2},a_{3}$ liegt dann $u_{12}$ in dem von $a_{1},a_{2}$ erzeugten Ideal in $K[V]^{G}$, d.h. es gibt $f_{1},f_{2} \in K[V]^{G}$ mit $a_{1}b_{2}-a_{2}b_{1}=f_{1}a_{1}+f_{2}a_{2}$. Weiter sind wegen der Regularit"at $a_{1},a_{2}$ teilerfremde Polynome\footnote{Ist  $d \in K[V]$ ein gemeinsamer Teiler, $a_{1}=dh$ so folgt $a_{2}h \in (a_{1})$, also gilt mit der Regularit"at $h \in (a_{1})$ , d.h. $d$ ist Einheit. Dabei beachte man, dass $a_{1},a_{2}$ ein phsop in dem CM Ring K[V] bilden, also dort sogar eine regul"are Sequenz.}, und es folgt, dass $a_{1}$ Teiler von $f_{2}+b_{1}$ ist, also $f_{2}+b_{1}=a_{1} \cdot h$ mit $h \in K[V]$. Nun ist
\[
a_{1} \cdot (\sigma -1)h=(\sigma -1)(a_{1}h)=(\sigma -1)(f_{2}+b_{1})=(\sigma -1)b_{1}=a_{1}g_{\sigma} \qquad \forall \sigma \in G,
\]
also $g_{\sigma}=(\sigma -1)h$. Damit ist die Restklasse von $g$ in $H^{1}(G,K[V])$ gleich 0, was die Injektivit"at beweist. \hfill $\Box$\\
\\
Wir k"onnen nun den zentralen Satz formulieren.\\
\\
\noindent {\bf Satz 6. } \textit{ Sei $G$ eine reduktive, nicht linear reduktive Gruppe.  Dann gibt es einen $G$-Modul $V$, so dass $K[V]^{G}$ nicht Cohen-Macaulay ist.\\
Ist genauer $0 \rightarrow U \rightarrow \tilde{U} \rightarrow K \rightarrow 0$ eine nicht zerfallende kurze exakte Sequenz von $G$-Moduln (die wegen Proposition 2 existiert) und $W$ ein $G$-Modul wie in Proposition 3, z.B. $W=\tilde{U}^{*}$, so kann man
\[
\begin{array}{rcl}
V & = & U^{*} \oplus W^{*} \oplus W^{*} \oplus W^{*}\\
(\textrm{bzw. } V & = & U^{*} \oplus \tilde{U} \oplus \tilde{U} \oplus \tilde{U})\\
\end{array}
\]
w"ahlen.}\\
\\
\noindent \textit{Beweis. } Sei also $0 \ne g \in H^{1}(G,U)$ und $W$ ein $G$-Modul, der ein $0 \ne w \in W^{G}$ mit $w \otimes g =0$ in $H^{1}(G, W \otimes U)$ enth"alt. Da $U^{**} \cong U, W^{**} \cong W$ haben wir also bis auf Isomorphie
\[
R:=K[V] \cong S(U \oplus W \oplus W \oplus W)
\]
(die symmetrische Algebra in $U \oplus W \oplus W \oplus W$, also  ``Polynome in Basis\-vek\-toren aus $U \oplus W \oplus W \oplus W$''). Sei $a_{i}$ die Kopie von $w$  im $i$-ten Summanden von $W$ in der direkten Summe ($i=1,2,3$). Dann liegen $a_{1},a_{2},a_{3}$ in $K[V]^{G}$ und bilden ein phsop in $K[V]$. Da $U$ ein direkter Summand von $R$ ist, haben wir eine Einbettung $H^{1}(G,U) \hookrightarrow H^{1}(G,R)$. Aufgrund der drei Kopien von $W^{*}$ in $V$ haben wir ebenfalls drei Einbettungen $H^{1}(G, W \otimes U) \hookrightarrow H^{1}(G,R)$. Das Bild von $w\otimes g$ unter diesen Einbettungen ist $a_{i} \cdot g$ ($i=1,2,3$). Da jedoch $w \otimes g=0$, wird das Bild von $g$ in $H^{1}(G,R)$ unter der durch Multiplikation mit $(a_{1},a_{2},a_{3})$ induzierten Abbildung $H^{1}(G,R) \rightarrow H^{1}(G,R^3)$ auf $0$ abgebildet. Wegen Proposition 5 kann daher $K[V]^{G}$ nicht Cohen-Macaulay sein. \hfill $\Box$

\section{Konstruktion nicht zerfallender Sequenzen}
Wir kommen nun zu dem Hauptresultat dieser Arbeit, n"amlich der Konstruk\-tion nicht zerfallender kurzer exakter Sequenzen f"ur die (reduktiven) Gruppen $GL_{n}(K)$ und  $SL_{n}(K)$ f"ur beliebiges $\textrm{char }K=p$.

Um die Bedeutung der Voraussetzung des n"achsten Satzes zu sehen, berechnen wir f"ur $\textrm{char } K=2$ die Menge aller $A \in GL_{2}(K)$ mit $A^{T}A=I_{2}$ (in positiver Charakteristik ist dies nicht $O_{2}(K)$) :\\
Aus $
\left(\begin{array}{cc}
a&b\\
c&d\\
\end{array}\right) 
\left(\begin{array}{cc}
a&c\\
b&d\\
\end{array}\right) 
=
\left(\begin{array}{cc}
1&0\\
0&1\\
\end{array}\right)
$ folgt $a^{2}+b^{2}=(a+b)^{2}=1$, also $b=a+1$, und analog $c=a+1$, sowie zuletzt noch $d=b+1=a$. Also ist die gesuchte Menge gleich
\[
\left\{ \left(\begin{array}{cc}
a&a+1\\
a+1&a\\
\end{array}\right) : \quad a\in K \right\} \textrm{       f"ur char }K=2 .
\]\\
\\
\noindent {\bf Satz 7. } \textit{ Sei $K$ ein (algebraisch abgeschlossener) K"orper mit $\textrm{char }K=p>0$, $n \ge 2$ und $G$ eine Untergruppe von $\textrm{GL}_{n}(K)$ mit\\
(a)  Falls $p=2$:  
$
\left(
\begin{array}{ccc}
 a & a+1 &\\
a+1&a&\\
&&I_{n-2}
\end{array} \right)
\in G
$
      f"ur wenigstens drei verschiedene Werte von $a\in K$ ($a=1$ ist stets ein solcher).\\
(b) Falls $p \ge 3$: 
$
\left(
\begin{array}{ccc}
1  & 1 &\\
0&1&\\
&&I_{n-2}
\end{array} \right)
\in G.
$\\
Sei weiter $V$ der Raum aller homogenen Polynome vom Grad $p$ in den Variablen $x_{1},...,x_{n}$. Ist $(a_{ij}) \in G$, so ist durch $(a_{ij}) \cdot x_{j} = \sum_{i=1}^{n}a_{ij}x_{i}$ die kanonische Operation von $G$ auf $V$ gegeben, wodurch $V$ zu einem $G$-Modul wird. Sei $W:=<x_{1}^{p},...,x_{n}^{p}> \le V$ und $U:=\{ f \in \textrm{Hom}_{K}(V,W): f|_{W}=0\}$. Da $W$ aufgrund des Frobenius-Homomorphismus ein $G$-Untermodul von $V$ ist, ist $U$ ein $G$-Modul. Sei ferner $\iota \in \textrm{Hom}_{K}(V,W)$ gegeben durch $\iota |_{W}=\textrm{id}_{W}$ und $\iota$ gleich $0$ auf allen Monomen, die nicht in $W$ liegen. Sei $\tilde{U}:=U \oplus K \iota$ und $\pi: \tilde{U} \rightarrow K$ gegeben durch $\pi(u + \lambda \cdot \iota):=\lambda \quad \textrm{f"ur } u \in U, \lambda \in K$. Dann ist durch
\[
0 \rightarrow U  \rightarrow  \tilde{U} \stackrel{\pi}{\rightarrow} K \rightarrow 0
\]
eine kurze exakte Sequenz gegeben, die nicht zerf"allt.}\\
\\
\noindent \textit{Beweis. } Wir verwenden f"ur $V$ eine monomiale Basis $\mathcal{B}$, wobei wir die Reihen\-fol\-ge der ersten $n+1$ bzw. $n+2$ Monome in den F"allen (a) bzw. (b) vorge\-ben, und zwar\\
im Fall (a):
\[
\mathcal{B}=\{x_{1}^{2},...,x_{n}^{2},x_{1}x_{2},...\}
\]\\
im Fall (b):
\[
\mathcal{B}=\{x_{1}^{p},...,x_{n}^{p},x_{1}^{p-1}x_{2},x_{1}^{p-2}x_{2}^{2},...\}
\]\\
Als Basis von $W$ dienen die ersten $n$ Eintr"age von $\mathcal{B}$. Sei $N:=|\mathcal{B}|= {n+p-1 \choose p}$. Wir bezeichnen mit $f_{p}:\textrm{GL}_{n}(K) \rightarrow \textrm{GL}_{n}(K)$ den koeffizientenweisen Frobenius-Homomorphismus, also $f_{p}(a_{ij})=(a_{ij}^{p})$. Ist $A_{\sigma} \in K^{N \times N}$ die Dar\-stell\-ungs\-matrix von $\sigma \in G$ bzgl. der Basis $\mathcal{B}$, so hat diese die Form
\[
A_{\sigma}=\left( \begin{array}{cc}
f_{p}(\sigma) & *\\
0&*\\
\end{array} \right).
\]
Weiter haben wir bzgl. der Basis $\mathcal{B}$
\[
U \cong \left\{
\left( \begin{array}{cc}
0_{n \times n} & B\\
\end{array} \right) \in K^{n \times N} \textrm{  mit  } B \in K^{n \times (N-n)} \right\},
\]
und f"ur $U \ni f \cong   
\left( \begin{array}{cc}
0_{n \times n} & B\\
\end{array} \right)$ bzgl. $\mathcal{B}$ haben wir die Operation gegeben durch
\[
\sigma \cdot f=\sigma \circ f \circ \sigma^{-1} \cong f_{p}(\sigma) \cdot
\left( \begin{array}{cc}
0_{n \times n} & B\\
\end{array} \right)
\cdot A_{\sigma^{-1}}.
\]
Die Darstellungsmatrix von $\iota$ ist gegeben durch
\[
\iota \cong
\left( \begin{array}{cc}
I_{n} & 0\\
\end{array} \right)=:J \in K^{n \times N}.
\]
Da
\[
\sigma \cdot \iota=\sigma \circ \iota \circ \sigma^{-1} \cong f_{p}(\sigma)
\left( \begin{array}{cc}
I_{n} & 0\\
\end{array} \right) 
\underbrace{
\left( \begin{array}{cc}
f_{p}(\sigma^{-1}) & *\\
0&*\\
\end{array} \right)}_{A_{\sigma^{-1}}} 
=
\left( \begin{array}{cc}
I_{n} & *\\
\end{array} \right)
\]
ist $\sigma \iota - \iota \in U$. Damit ist $\pi$ wohldefiniert, und mit $g_{\sigma}:=(\sigma -1)\iota$ ist $g \in Z^{1}(G,U)$. Wir m"ussen zeigen, dass $g \not\in B^{1}(G,U)$. Wir nehmen das Gegenteil an, also die Existenz eines
\[
U \ni u \cong Z =
\left( \begin{array}{cc}
0_{n \times n} & \hat{Z}\\
\end{array} \right) \in K^{n \times N}, \qquad \hat{Z}=(z_{ij}) \in K^{n \times (N-n)}
\]
mit $g_{\sigma}=(\sigma -1) \iota \stackrel{!}{=}(\sigma -1)u \quad \textrm{f"ur alle } \sigma \in G$ bzw.
\begin{equation} \label{zerf}
f_{p}(\sigma)J A_{\sigma^{-1}}-J \stackrel{!}{=}f_{p}(\sigma)Z A_{\sigma^{-1}}-Z \quad \forall \sigma \in G.
\end{equation}
Diese letzte Gleichung f"uhren wir nun in beiden F"allen zum Widerspruch.\\

\noindent (a) Mit
\[
\sigma^{-1}=
\left(
\begin{array}{ccc}
 a & a+1 &\\
a+1&a&\\
&&I_{n-2}
\end{array} \right)
\]
berechnen wir die $(n+1)$te Spalte von $A_{\sigma^{-1}}$:
\[
\begin{array}{rcl}
\sigma^{-1} \cdot x_{1}x_{2}&=&(ax_{1}+(a+1)x_{2})((a+1)x_{1}+ax_{2})\\
&=&(a^{2}+a)x_{1}^{2}+(a^{2}+a)x_{2}^{2}+x_{1}x_{2}\\
&\cong&(a^{2}+a,a^{2}+a,0_{n-2},1,0)^{T}.\\
\end{array}
\]
Damit vergleichen wir nun auf beiden Seiten von (\ref{zerf}) den Eintrag in der ersten Zeile und $(n+1)$ten Spalte:\\
Links:
\[
(a^{2},a^{2}+1,0_{n-2})
\left( \begin{array}{cc}
I_{n} & 0_{n \times (N-n)}\\
\end{array} \right)
\left( \begin{array}{c}
a^{2}+a\\
a^{2}+a\\
0_{n-2}\\
1\\
0\\
\end{array} \right) =a^{2}+a
\]
Rechts:
\[
\begin{array}{rcl}
(a^{2},a^{2}+1,0_{n-2})
\left( \begin{array}{cc}
0_{n \times n} & \hat{Z}\\
\end{array} \right)
\left( \begin{array}{c}
a^{2}+a\\
a^{2}+a\\
0_{n-2}\\
1\\
0\\
\end{array} \right) -z_{11}&=&a^{2}z_{11}+(a^{2}+1)z_{21}-z_{11}\\
&=&(a^{2}+1)(z_{11}+z_{21})\\
\end{array}
\]
Setzen wir $c:=z_{11}+z_{21}$, so sehen wir, dass $ca^2+c=a^{2}+a$ bzw.
\[
(c+1)a^{2}+a+c=0
\]
f"ur wenigstens drei verschiedene Werte $a \in K$ erf"ullt sein muss. Dies ist ein Widerspruch.\\

\noindent (b) Wir betrachten
\[
\sigma=
\left(
\begin{array}{ccc}
1  & 1 &\\
0&1&\\
&&I_{n-2}
\end{array} \right),
\sigma^{-1}=
\left(
\begin{array}{ccc}
1  & -1 &\\
0&1&\\
&&I_{n-2}
\end{array} \right)
\in G
\]
und berechnen die $(n+1)$te und $(n+2)$te Spalte von $A_{\sigma^{-1}}$:\\
$(n+1)$te Spalte:
\[
\begin{array}{rcl}
\sigma^{-1} \cdot x_{1}^{p-1}x_{2} &=& x_{1}^{p-1}(-x_{1}+x_{2})\\
&=&-x_{1}^{p}+x_{1}^{p-1}x_{2}\\
&\cong &(-1,0_{n-1},1,0)^{T}
\end{array}
\]
$(n+2)$te Spalte:
\[
\begin{array}{rcl}
\sigma^{-1} \cdot x_{1}^{p-2}x_{2}^{2} &=& x_{1}^{p-2}(x_{1}^{2}-2x_{1}x_{2}+x_{2}^{2})\\
&=&x_{1}^{p}-2x_{1}^{p-1}x_{2}+x_{1}^{p-2}x_{2}^{2}\\
&\cong &(1,0_{n-1},-2,1,0)^{T}
\end{array}
\]
Wir vergleichen nun wieder beide Seiten von (\ref{zerf}):\\
(i) \underline{erste Zeile, $(n+1)$te Spalte}\\
Links:
\[
\left( \begin{array}{ccc} 1 & 1 & 0_{n-2}\\ \end{array} \right)
\left( \begin{array}{cc} I_{n} & 0_{n \times (N-n)}\\ \end{array} \right)
\left( \begin{array}{c}
-1\\
0_{n-1}\\
1\\
0\\
\end{array} \right)=-1
\]
Rechts:
\[
\left( \begin{array}{ccc} 1 & 1 & 0_{n-2}\\ \end{array} \right)
\left( \begin{array}{cc} 0_{n \times n} & \hat{Z}\\ \end{array} \right)
\left( \begin{array}{c}
-1\\
0_{n-1}\\
1\\
0\\
\end{array} \right)-z_{11}=z_{11}+z_{21}-z_{11}=z_{21}
\]
Da Gleichheit gelten soll haben wir
\begin{equation} \label{wid1}
z_{21}=-1
\end{equation}

\noindent (ii) \underline{zweite Zeile, $(n+2)$te Spalte}\\
Links:
\[
\left( \begin{array}{ccc} 0 & 1 & 0_{n-2}\\ \end{array} \right)
\left( \begin{array}{cc} I_{n} & 0_{n \times (N-n)}\\ \end{array} \right)
\left( \begin{array}{c}
1\\
0_{n-1}\\
-2\\
1\\
0\\
\end{array} \right)=0
\]
Rechts:
\[
\left( \begin{array}{ccc} 0 & 1 & 0_{n-2}\\ \end{array} \right)
\left( \begin{array}{cc} 0_{n \times n} & \hat{Z}\\ \end{array} \right)
\left( \begin{array}{c}
1\\
0_{n-1}\\
-2\\
1\\
0\\
\end{array} \right)-z_{22}=-2z_{21}+z_{22}-z_{22}=-2z_{21}
\]
Da $p \ge 3$ ist $2 \ne 0$, und der Vergleich beider Seiten liefert
\[
z_{21}=0,
\]
im Widerspruch zu (\ref{wid1}). \hfill$\Box$\\
\\
Mit obigen Bezeichnungen ist also wegen Satz 6 mit 
\[
X:=U^{*} \oplus \tilde{U} \oplus \tilde{U} \oplus \tilde{U}
\]
durch $K[X]^{G}$
ein Invariantenring ohne die Cohen-Macaulay Eigenschaft ge\-ge\-ben, wobei 
\[
\textrm{dim }X=4\cdot\textrm{dim }U+3=4 \cdot n(N-n)+3=4n \left( {n+p-1 \choose p} - n\right)+3.
\]
Insbesondere f"ur $n=2$ haben wir $\textrm{dim }X=8p-5$, also hat der kleinste so konstruierbare Modul $X$ die Dimension $11$ f"ur $p=2$.
Dabei erlauben die Voraussetzungen von Satz 7 f"ur $G$ insbesondere jede reduktive Zwischen\-grup\-pe $\textrm{SL}_{n}(K) \le G \le \textrm{GL}_{n}(K)$. F"ur $p \ge 3$ ist au\ss erdem die von der Matrix im Fall (b) erzeugte Untergruppe isomorph zu der (reduktiven) endlichen Gruppe $\mathbf{Z}_{p}$, und noch weitere endliche Untergruppen von $\textrm{GL}_{n}(K)$ erf"ullen die Voraussetzungen von (a) bzw. (b), z.B. hat man im Fall (b) mit $n = 4$ und dem Primk"orper $\mathbf{F}_{p}$ mit
\[
G=\left\{
\left(\begin{array}{cccc}
1 &a\\
0&1\\
&&1&b\\
&&0&1\\
\end{array}\right): a,b \in \mathbf{F}_{p}
\right\}
\]
eine Gruppe, die isomorph zu $\mathbf{Z}_{p} \times \mathbf{Z}_{p}$ ist.\\
Sei allgemeiner $U$ ein $k$-dimensionaler Untervektorraum des un\-end\-lich di\-men\-sion\-al\-en $\mathbf{F}_{p}$-Vektorraums $K$, wobei $\textrm{char }K=p \ge 2$. Setzt man
\[
A(a):=\left\{
\begin{array}{cl}
\left(\begin{array}{cc}
a+1 &a\\
a&a+1\\
\end{array}\right) &\textrm{f"ur }p =2 \\
\left(\begin{array}{cc}
1 &a\\
0&1\\
\end{array}\right) & \textrm{f"ur }p \ge 3 \\
\end{array}
\right.
\]
so gilt $A(a)\cdot A(b)=A(a+b)$. Daher ist $G:=\left\{A(a): a \in U\right\} \cong \mathbf{Z}_{p}^{k}$, wie man anhand einer Basisdarstellung der Elemente aus $U$ sieht. Damit die Voraussetzungen von Satz 7 erf"ullt sind, muss $k \ge 2$ f"ur $p=2$ bzw. $1 \in U$ f"ur $p \ge 3$ gelten.\\
\\
Zum Schlu\ss{ }sei noch auf einen anderen Ansatz zur Konstruktion einer nicht zerfallenden, kurzen exakten Sequenz hingewiesen, der jedoch nur im Falle $p=n=2$ zum Ziel f"uhrte. Der Einfachheit halber sei $G=\textrm{SL}_{2}(K)$. Man betrachtet nun die kanonische Operation von $G$ auf den homogenen Polynomen vom Grad $2$ in zwei Variablen $x$ und $y$. Setzt man $U:=<x^{2},y^{2}>$, $\tilde{U}:=<x^{2},y^{2},xy>$ und $\pi: \tilde{U} \rightarrow K, \quad \pi(ax^{2}+by^{2}+cxy):=c$, so ist $0 \rightarrow U  \rightarrow  \tilde{U} \rightarrow K \rightarrow 0$ eine nicht zerfallende kurze exakte Sequenz. Durch Tensorierung mit $\textrm{det}^{-1}$ kann man ein Beispiel f"ur $\textrm{GL}_{2}(K)$ angeben. Wie gesagt l"asst sich dieses Beispiel jedoch nicht auf $p \ge 3$ bzw. $n \ge 3$ verallgemeinern. Abgesehen davon rechnet man leicht nach, dass die so konstruierte Sequenz zu der von Satz 7 mit $n=p=2$ "aquivalent ist.

\end{document}